# EMBEDDED SPHERES AND 4–MANIFOLDS WITH SPIN COVERINGS

CHRISTIAN BOHR

ABSTRACT. A strategy for constructing an embedded sphere in a 4–manifold realizing a given homology class which has been successfully applied in the past is to represent the class as a first step stably by an embedded sphere, i.e. after adding products of 2–spheres, and to move that sphere back into the original manifold. In this paper, we study under what conditions the first step of this approach can be carried out if the 4–manifold at hand is not simply connected. One of our main results is that there are – apart from the well known Arf invariant – additional bordism theoretical obstructions to stably representing homology classes by embedded spheres.

## 1. INTRODUCTION AND SUMMARY OF RESULTS

It has always been one of the most challenging problems in 4–dimensional topology to determine which homology classes in 4–manifolds can be represented by embedded spheres. Although considerable progress has been made during the last decades, this question is in general still open. An important contribution was made by R. Lee and D. Wilczyński in a series of papers [7] [8] in which they gave necessary and sufficient conditions for a homology class $\alpha \in H_2(X;\mathbb{Z})$, where $X$ is a simply connected 4–manifold, to be representable by a locally flat topological embedding of a 2–sphere. To construct such an embedded sphere, they first showed that the class can be stably represented by an embedded sphere, i.e. after adding sufficiently many copies of $S^2 \times S^2$, and then proved that this sphere can be moved into $X$ itself by a surgery argument. If one aims at extensions of these results to non–simply connected 4–manifolds using a similar approach, one therefore has to study the question whether a given homology class can be stably represented by an embedded sphere without assuming that $X$ be simply connected, which is the main objective of this paper.

To formulate this question more precisely, assume that $X$ is a 4–manifold and that $\xi \in H_2(X;\mathbb{Z})$ is some homology class. For every non–negative integer $k$ we can then consider the 4–manifold

$$X_k = X \# k(S^2 \times S^2)$$

obtained from $X$ by adding $k$ copies of $S^2 \times S^2$. In a natural way, $H_2(X;\mathbb{Z})$ is a subgroup of $H_2(X_k;\mathbb{Z})$ and therefore we can think of $\xi$ as a homology class in $X_k$. We will say that the class $\xi$ can be stably represented by an embedded sphere if there is, for some $k$, an embedded sphere in $X_k$ representing $\xi$. For the sake of simplicity we will assume that $X$ is smooth and that all the embeddings are smooth embeddings, although there are topological versions of most of our results.

1991 *Mathematics Subject Classification.* 57M99;55N22.

The author has been supported by a fellowship from the *Deutsche Forschungsgemeinschaft*.

Of course a class which can be stably represented by an embedded sphere is spherical, i.e. it is in the image of the natural map from the second homotopy group to the second homology group, therefore we will restrict our attention to those classes. It was understood early one that one has to distinguish between those classes $\alpha$ which are characteristic, i.e. the Poincaré dual $PD(\alpha)$ is a lift of the second Stiefel–Whitney class, and those classes for which this is not the case. In the early sixties, M. Kervaire and J. Milnor discovered in [5] that Rokhlin's Theorem on the signatures of smooth 4–dimensional spin manifolds implies an obstruction to representing characteristic homology classes by embedded spheres As shown later by M. Freedman and R. Kirby [3] and independently by L. Taylor, the vanishing of this obstruction is not only necessary but also sufficient, provided that the 4–manifold $X$ is simply connected.

**Theorem 1** (Freedman, Kirby, Taylor). *Suppose that $X$ is a simply connected 4–manifold and that $\xi \in H_2(X;\mathbb{Z})$ is a characteristic class. Then $\xi$ can be stably represented by an embedded sphere if and only if $\xi \cdot \xi \equiv \sigma(X) \mod 16$.*

An obvious approach to constructing a sphere representing $\xi$ is to start with an embedded surface $F$ of arbitrary genus and to remove handles to reduce the genus. Theorem 1 is proved by showing that this approach works if and only if a certain obstruction, the so called Arf invariant of $F$, vanishes and that this Arf invariant which takes values in $\mathbb{Z}_2$ is given by the formula

$$\operatorname{Arf}(F) = \frac{1}{8}(\xi \cdot \xi - \sigma(X)) \mod 2.$$

This formula provides a relation between Rokhlin's Theorem, which can be recovered from Theorem 1 by applying it the the trivial homology class if $X$ is spin, and geometric properties of characteristic surfaces in 4–manifolds. Finding similar relations between bordism theory and a larger class of embedded surfaces is the main intention of this paper.

It is natural to ask for extensions of Theorem 1 to the case of non–simply connected 4–manifolds. In order to be able to perform surgery on an embedded surface representing $\xi$ aiming to reduce its genus, it is desirable to have an embedding $F \to X_k$ for some $k$ such that the induced map $\pi_1(F) \to \pi_1(X_k)$ is trivial. Such an embedding is usually called a $\pi_1$–null embedding, and if it exists we will say that the class $\xi$ can be stably represented by a $\pi_1$–null embedding. As shown in [1], this is always the case if $\xi$ is spherical and characteristic, which can be used to extend Theorem 1 to the non–simply connected case.

**Theorem 2** ([1]). *Suppose that $X$ is a 4–manifold and that $\xi \in H_2(X;\mathbb{Z})$ is a spherical characteristic class. Then $\xi$ can be stably represented by a $\pi_1$–null embedded surface. Moreover that surface can be chosen to be a sphere if and only if $\xi \cdot \xi \equiv \sigma(X) \mod 16$.*

These results give a rather satisfactory answer to our question if the homology class at hand is characteristic. In the presence of a non–trivial fundamental group, there are other distinguished homology classes which deserve special attention, namely those homology classes whose Poincaré dual becomes characteristic if pulled back to the universal covering. We will call these classes spherically characteristic.

**Definition 1.** Suppose that $X$ is a 4–manifold with universal covering $\pi\colon \widetilde{X} \to X$. A homology class $\alpha \in H_2(X;\mathbb{Z})$ is called **spherically characteristic** if

$$\pi^* PD(\alpha) = w_2(\widetilde{X}),$$

where $PD$ denotes the Poincaré duality map. A class which is not spherically characteristic will be called **ordinary**.

Similarly to the correspondence between characteristic classes and spin manifolds indicated above, these classes correspond to 4–manifolds whose universal covering is spin. In fact, the universal covering of a 4–manifold $X$ is spin if and only if the trivial homology class is spherically characteristic.

As to the case of ordinary classes, it is not difficult to prove that these classes can always be stably represented by embedded spheres. This fact was well known in the simply connected case, in the general case it easily follows from the arguments given in [2].

**Theorem 3** ([2]). *Every ordinary spherical homology class can be stably represented by an embedded sphere.*

However this is not true for classes which are not characteristic, but spherically characteristic. In fact it turns out that some of these classes cannot even be stably represented by a $\pi_1$–null embedding. This is a consequence of a more general result (Theorem 7) stated in Section 4, which describes a bordism theoretical obstruction whose vanishing is a necessary and sufficient condition for the existence of a stable $\pi_1$–null embedding representing a given spherically characteristic and spherical homology class.

**Theorem 4.** *There is a 4–manifold $X$ and a spherically characteristic spherical class $\xi \in H_2(X;\mathbb{Z})$ which cannot be stably represented by a $\pi_1$–null embedding. The fundamental group of $X$ can be chosen to be $\mathbb{Z} \oplus \mathbb{Z}_2$.*

Now suppose that we can actually find a $\pi_1$–null embedded surface $F$ representing a spherically characteristic class, we will see later that there are many fundamental groups for which this is always the case (Theorem 8). Then we can again try to reduce the genus of that surface to obtain an embedded sphere. As for characteristic surfaces, one can show that this works if and only if the Arf invariant $\text{Arf}(F)$ (whose definition has to be adapted slightly to the more general setting) vanishes. We have seen above that Rokhlin's Theorem implies that certain homology classes cannot be represented by surfaces with vanishing Arf invariant. However Rokhlin's Theorem fails for the larger class of 4–manifolds whose universal covering is spin, as the following result due to P. Teichner [11] shows.

**Theorem 5** (Teichner). *Every integer can be realized as the signature of a 4–manifold whose universal covering is spin.*

In contrast to the characteristic case, there is therefore no obvious reason why it should not be possible to modify a given $\pi_1$–null embedding in such a way that the Arf invariant becomes zero, but the homology class does not change. Thus one expects that once we have found a $\pi_1$–null embedded representative, there should in general be no further obstructions to representing that homology class by an embedded sphere. In fact this turns out to be true, at least if the homology class at hand is in a certain sense primitive.

**Theorem 6.** *Suppose that $X$ is a 4–manifold and that $F \hookrightarrow X$ is a $\pi_1$–null embedding which is spherically characteristic, but not characteristic. Assume further that there is a spherical class $\omega \in H_2(X; \mathbb{Z})$ such that $\omega \cdot F = 1$. Then there is, for some non–negative $k$, a $\pi_1$–null embedding*

$$F' \longrightarrow X \# k(S^2 \times S^2)$$

*such that $[F'] = [F]$, but $\mathrm{Arf}(F') \neq \mathrm{Arf}(F')$. In particular the class $[F]$ can be stably represented by an embedded sphere.*

We believe that the presence of the class $\omega$ is not necessary, i.e. that the statement of the Theorem also holds without assuming the existence of such a class, although our methods do not work without this assumption. The main idea of the proof is to use the rich supply of 4–manifolds whose universal covering is spin guaranteed by Teichner's results to construct embedded surfaces with prescribed Arf invariant. This interplay between spherically characteristic surfaces and 4–manifolds with spin coverings generalizes the well known relations between characteristic surfaces and spin 4–manifolds.

## 2. An Arf invariant for spherically characteristic surfaces

In this section, we will describe an extension of the Arf invariant of characteristic surfaces in 4–manifolds to surfaces representing spherically characteristic homology classes. Descriptions of the Arf invariant in the characteristic case can be found in [3], [9] and [4], where it is used to give a geometric proof of Rokhlin's Theorem. Our outline closely follows the presentation in [3], therefore we will focus on the points where modifications are needed and refer the reader to [3] for technical details. Unless stated otherwise, all manifolds in this and the following sections will be assumed to be closed, connected, oriented and smooth.

Of course an embedded surface $F$ in a 4–manifold $X$ is called spherically characteristic if its homology class is spherically characteristic in the sense of Definition 1.

**Definition 2.** A **spherically characteristic 2–ad** is a quadruple $(A, B, F, X)$, where $X$ is a 4–manifold, $F \subset X$ a spherically characteristic surface, $A \subset F$ is a closed submanifold of dimension 1 and $B \subset X$ is a connected surface which meets $F$ normally along $A = \partial B$ such that the interior of $B$ intersects $F$ transversely in $F \setminus A$ and such that the map $\pi_1(B) \longrightarrow \pi_1(X)$ is trivial.

**Definition 3.** We will say that two spherically characteristic 2–ads $(A, B, F, X)$ and $(A', B', F', X')$ are **equivalent** if there is an oriented connected bordism $(\bar{X}, \bar{F})$ with boundary $\partial(\bar{X}, \bar{F}) = (X, F) - (X', F')$ such that $\bar{F} \subset \bar{X}$ is spherically characteristic, together with a surface $\bar{A} \subset \bar{F}$, not necessarily orientable, with $\partial \bar{A} = A \cup A'$ such that the homomorphism $\pi_1(\bar{A}_i) \to \pi_1(\bar{X})$ is trivial for every component $\bar{A}_i$ of $\bar{A}$. Here $\bar{F}$ is said to be spherically characteristic if the pullback of the cohomology class $PD([\bar{F}]) \in H^2(\bar{X}; \mathbb{Z})$ to the universal covering of $\bar{X}$ is characteristic.

Assume that we are given a spherically characteristic 2–ad $(A, B, F, X)$. Let $N_B^X$ denote the normal bundle of $B$ in $X$. Pick a section without zeros $V$ of the normal bundle $N_A^F \to A$ of $A$ in $F$. Then $V$ defines a section of the restriction $N_B^X|A$. We have a well–defined obstruction $w_2(N_B^X; V) \in \mathbb{Z}_2$ to extending $V$ to a section of $N_B^X$ which is nowhere zero.

**Definition 4.** Let us a define an invariant $q$ of a spherically characteristic 2–ad $(A, B, F, X)$ by setting
$$q(A, B, F, X) = w_2(N_B^X; V) + \#(int(B) \cap F) \mod 2.$$

We will now show that the value of $q$ only depends on the homology class of $A$. As in [3], the key step is the following observation.

**Lemma 1.** *If $(A, B, F, X)$ and $(A', B', F', X')$ are equivalent, then*
$$q(A, B, F, X) = q(A', B', F', X').$$

*Proof.* By the very definition of equivalence, there is a connected bordism $(\bar{X}, \bar{F})$ between the pairs $(X, F)$ and $(X', F')$ such that $\bar{F}$ is spherically characteristic, and a surface $\bar{A} \subset \bar{F}$ with $\partial \bar{A} = A - A'$ such that every component of $\bar{A}$ is a $\pi_1$–null embedding. For the sake of simplicity we will restrict ourselves to the case that $\bar{A}$ is connected, it is clear how the proof has to be modified in the general case.

Let $Y = B \cup \bar{A} \cup B'$ denote the surface which is obtained by gluing together $B, B'$ and $\bar{A}$ and smoothing the angle. By the van Kampen–Theorem, the map $\pi_1(Y) \longrightarrow \pi_1(\bar{X})$ is trivial, as this is the case for every single piece. We can now prove our claim using exactly the same arguments as in [3] once we show that
$$\bar{F} \cdot Y = \langle w_2(\bar{X}), [Y] \rangle \in \mathbb{Z}_2.$$
But this equality holds true because $\bar{F}$ is spherically characteristic and the embedding $Y \to \bar{X}$ admits a lift to the universal covering of $\bar{X}$. $\square$

Now suppose we are given a 4–manifold $X$ and a spherically characteristic $\pi_1$–null embedded surface $F \subset X$. Let $A \subset F$ be a 1–dimensional closed submanifold. Since, by assumption, the homotopy class of $A$ in $X$ is zero, there exists a $\pi_1$–null embedded surface $B \subset X$ with boundary $\partial B = A$, we can even assume that $B$ is an embedded disk. After moving $B$ into general position, $(A, B, F, X)$ will be a spherically characteristic 2–ad. Let $q(A) = q(A, B, F, X)$.

**Lemma 2.** *The value of $q(A)$ is independent of the choice of the surface $B$ and only depends on the $\mathbb{Z}_2$–homology class of $A$.*

*Proof.* Suppose that $A'$ is another representative of $[A] \in H_1(F; \mathbb{Z}_2)$ and $B' \subset X$ is a $\pi_1$–null embedded surface with boundary $A'$. Then there exists a – possibly non–orientable – surface
$$\bar{A} \subset \bar{F} = F \times [0, 1]$$
joining $A$ and $A'$. Since the embedding $\bar{F} \to \bar{X} = X \times [0, 1]$ is a $\pi_1$–null embedding, the same is true for every component of $\bar{A}$. Therefore $(A, B, F, X)$ and $(A', B', F, X)$ are equivalent and the claim follows. $\square$

Thus $q$ defines a map $q \colon H_1(X; \mathbb{Z}_2) \to \mathbb{Z}$, which, by the same arguments as in [3], is a quadratic form.

**Definition 5.** Let $X$ be a 4–manifold and $F \subset X$ a spherically characteristic $\pi_1$–null embedded surface. Define the **Arf invariant** $\mathrm{Arf}(F) \in \mathbb{Z}_2$ to be the Arf invariant of the quadratic form $q$ constructed above.

It is easy to check that the Arf invariant is – as in the characteristic case – additive under connected sum. More precisely let us assume that $F$ and $F'$ are $\pi_1$–null embedded spherically characteristic surfaces in 4–manifolds $X$ and $X'$.

Then we have a spherically characteristic $\pi_1$–null embedding $F \# F' \subset X \# X'$ and $\mathrm{Arf}(F \# F') = \mathrm{Arf}(F) + \mathrm{Arf}(F')$. Moreover the Arf invariant of a spherically characteristic $\pi_1$–null embedded surface is again the only obstruction to performing surgery to obtain an embedded sphere representing the same homology class, i.e. we have the following lemma which can easily be proved by adapting the arguments given in [3], see also [1] for a slightly different proof.

**Lemma 3.** *Let $F \subset X$ be a surface of positive genus in a 4–manifold such that $F$ is spherically characteristic, $\pi_1(F) \to \pi_1(X)$ is trivial and $\mathrm{Arf}(F) = 0$. Then there is a non–negative integer $k$ and an embedded sphere in $X \# k(S^2 \times S^2)$ representing the class $[F]$.*

Adapting the arguments used in [12], one can also see that starting with a surface of genus $g$, it suffices to add $g$ copies of $S^2 \times S^2$ in order to obtain an embedded sphere.

It is worth mentioning that, in contrast to the characteristic case, the Arf invariant of a spherically characteristic surface is in general not invariant under spherically characteristic bordisms. In fact, as a consequence of Theorem 6, the Arf invariant is not even an invariant of the homology class but depends on the way the surface is embedded in $X$.

## 3. Bordisms groups and spherically characteristic classes

In this section, we will introduce and investigate certain bordism groups and maps between them. These bordism groups are given by certain fibrations, generalizing the diagram of fibrations

$$\begin{array}{ccc} \mathrm{BSpin} & \longrightarrow & \mathrm{BSpin}^c \\ \downarrow & & \downarrow \\ \mathrm{BSO} & = & \mathrm{BSO} \end{array}$$

which are tailored such that spherically characteristic classes are characteristic classes of the corresponding normal structures. The motivating special case the reader should keep in mind is the case of characteristic classes which are exactly the Chern classes of Spin$^c$–structures. This will allow us to use Kreck's modified surgery theory [6] to construct stable diffeomorphisms with a prescribed behavior with respect to a given spherically characteristic class. Let us now start to make these ideas precise.

Suppose we are given a finitely presentable group $\Pi$ and a cohomology class $w \in H^2(\Pi; \mathbb{Z})$. Then we can form fibrations

$$K(\mathbb{Z}, 2) \longrightarrow B^c(\Pi, w) \xrightarrow{p} \mathrm{BSO} \times K(\Pi, 1)$$

and

$$K(\mathbb{Z}_2, 1) \longrightarrow B(\Pi, w) \xrightarrow{q} \mathrm{BSO} \times K(\Pi, 1)$$

by pulling back the path space fibrations over $K(\mathbb{Z}, 3)$ respectively $K(\mathbb{Z}_2, 2)$ using the maps given by the cohomology classes

$$\beta(w_2 + w) \in H^3(\mathrm{BSO} \times K(\Pi, 1); \mathbb{Z})$$

and

$$w_2 + w \in H^2(\mathrm{BSO} \times K(\Pi, 1); \mathbb{Z}_2).$$

Here $\beta$ denotes the Bockstein operator associated to the short exact sequence

$$0 \longrightarrow \mathbb{Z} \xrightarrow{\cdot 2} \mathbb{Z} \longrightarrow \mathbb{Z}_2 \longrightarrow 0.$$

Of course these fibrations are only determined up to fibre homotopy equivalence.

To justify our notation, we remark that the fibration $B(\Pi, w) \to \mathrm{BSO}$ given by composing $q$ and the projection to the first factor is fibre homotopy equivalent to the fibration which is called the "normal 1–type" of a 4–manifold with spin coverings in [6] and is also used in [11]. We will need some basic facts about $p$ and $q$ which are summarized in the following lemma.

**Lemma 4.**
1. *The restriction of $B^c(\Pi, w) \to \mathrm{BSO} \times K(\Pi, 1)$ to $\mathrm{BSO}$ is the usual fibration $\mathrm{BSpin}^c \to \mathrm{BSO}$. Consequently the fibre of the fibration*
   $$p' \colon B^c(\Pi, w) \longrightarrow K(\Pi, 1)$$
   *obtained by composing $p$ and the projection to the second factor is $\mathrm{BSpin}^c$. Moreover this fibration is orientable.*
2. *Similarly the restriction of $B(\Pi, w) \to \mathrm{BSO} \times K(\Pi, 1)$ to $\mathrm{BSO}$ is the usual fibration $\mathrm{BSpin} \to \mathrm{BSO}$. Hence the fibre of the fibration*
   $$q' \colon B(\Pi, w) \longrightarrow K(\Pi, 1)$$
   *obtained by composing $q$ and the projection to the second factor is $\mathrm{BSpin}$. Also this fibration is orientable.*
3. *The map $p'$ induces an isomorphism $\pi_1(B^c(\Pi, w)) = \Pi$. Similarly $q'$ induces an isomorphism $\pi_1(B(\Pi, w)) = \Pi$.*
4. *The homomorphism $H^2(B^c(\Pi, w); \mathbb{Z}) \to H^2(\mathrm{BSpin}^c; \mathbb{Z})$ induced by the restriction to the fibre of $p'$ is onto, and there is a lift $\hat{c}_1 \in H^2(B^c(\Pi, w); \mathbb{Z})$ of $p^*(w_2 + w)$ which restricts to $c_1 \in H^2(\mathrm{BSpin}^c; \mathbb{Z})$. For such a class $\hat{c}_1$, the homomorphism*
   $$\pi_2(B^c(\Pi, w)) \longrightarrow \mathbb{Z}$$
   *mapping $[f]$ to $\langle \hat{c}_1, f_*[S^2] \rangle$ is an isomorphism.*

*Proof.* The restriction of $B^c(\Pi, w) \to \mathrm{BSO} \times K(\Pi, 1)$ to $\mathrm{BSO}$ is the fibration determined by the cohomology class $W_3 = \beta(w_2)$ which is clearly $\mathrm{BSpin}^c \to \mathrm{BSO}$. Hence the fibre of $p'$ is $\mathrm{BSpin}^c$. The exact homotopy sequence of the fibration $p'$ shows that $\pi_1(B^c(\Pi, w)) = \Pi$, where the isomorphism is induced by $p'$.

To prove the remaining statements, let us describe another way of constructing $B^c(\Pi, w)$, which is also useful for other purposes. To this end let $f \colon E \to K(\Pi, 1)$ denote the pullback of the fibration

(*)  $$\mathrm{BSpin}^c \longrightarrow \mathrm{BSO} \longrightarrow K(\mathbb{Z}, 3)$$

via the map $K(\Pi, 1) \to K(\mathbb{Z}, 3)$ given by $\beta(w)$. Let $\hat{f} \colon E \to \mathrm{BSO}$ denote the canonical map. Then we have a map

$$g \colon E \longrightarrow \mathrm{BSO} \times K(\Pi, 1)$$

given by $f$ and $\hat{f}$ which we can turn into a fibration. Using standard tools like the exact homotopy sequence and the Leray–Serre spectral sequence, it is not difficult to check that the fibration obtained in this way is fibre homotopy equivalent to $p$. In particular the fibration $p'$ is the pullback of the fibration (*) and is therefore orientable.

Now let us turn to the proof of the corresponding statements for $q$. The restriction of that fibration to BSO is given by $w_2$ and is therefore fibre homotopy equivalent to BSpin $\to$ BSO. Consequently the fibre of $B(\Pi, w) \to K(\Pi, 1)$ is BSpin. The homotopy groups of $B(\Pi, w)$ are easily computed using the exact homotopy sequence, observe that the map $\pi_2(\mathrm{BSO} \times K) \to \pi_2(K(\mathbb{Z}_2, 2))$ induced by $w_2 + w$ is not zero. The fibration $q'$ can also be constructed by pulling back the fibration BSO $\to K(\mathbb{Z}_2, 2)$ determined by $w_2$ using the map $K(\Pi, 1) \to K(\mathbb{Z}_2, 2)$ given by $w$, which proves in particular that it is orientable.

As to statement 4, the description of $p'$ as the pullback of the fibration (*) (whose transgression maps $c_1$ to $2\iota$) shows that the transgression map of the fibration $p'$ is trivial, therefore the restriction map

$$\iota^* \colon H^2(B^c(\Pi, w); \mathbb{Z}) \longrightarrow H^2(\mathrm{BSpin}^c; \mathbb{Z})$$

is onto. By construction of the fibration $p$, $p^*\beta(w_2 + w) = 0$, hence $p^*(w_2 + w)$ has a lift $x \in H^2(B^c(\Pi, w); \mathbb{Z})$. The reduction modulo 2 of the restriction $\iota^* x$ to $\mathrm{BSpin}^c$ is $\iota^* p^* w_2$ which is also the reduction of $c_1$. As $\iota^*$ is onto, we can therefore add some multiple of 2 to $x$ to obtain a lift $\hat{c}_1$ with $\iota^* \hat{c}_1 = c_1$. To prove that the map $\pi_2(B^c(\Pi, w)) \to \mathbb{Z}$ given by $\hat{c}_1$ is an isomorphism, note that $\pi_2(\mathrm{BSpin}^c) \to \pi_2(B^c(\Pi, w))$ is an isomorphism, so we only have to prove that the map $\pi_2(\mathrm{BSpin}^c) \to \mathbb{Z}$ given by $\iota^* \hat{c}_1 = c_1$ is an isomorphism. But this follows from the Universal Coefficient Theorem and the Hurewicz Theorem. $\square$

**Definition 6.** An **admissible triple** is a triple $(\Pi, w, \hat{c}_1)$, where $\Pi$ is a finitely presentable group, $w \in H^2(\Pi; \mathbb{Z}_2)$ and

$$\hat{c}_1 \in H^2(B^c(\Pi, w); \mathbb{Z})$$

is a lift of $p^*(w_2 + w)$ restricting to $c_1 \in H^2(\mathrm{BSpin}^c; \mathbb{Z})$ as in Lemma 4.

**Lemma 5.** *Suppose that $(\Pi, w, \hat{c}_1)$ is an admissible triple. Then there is a fibration*

$$\pi \colon B(\Pi, w) \longrightarrow B^c(\Pi, w)$$

*with fibre $K(\mathbb{Z}, 1)$ such that the diagram*

$$\begin{array}{ccc} B(\Pi, w) & \xrightarrow{\pi} & B^c(\Pi, w) \\ {\scriptstyle q}\downarrow & & \downarrow{\scriptstyle p} \\ \mathrm{BSO} \times K(\Pi, 1) & =\!=\!= & \mathrm{BSO} \times K(\Pi, 1) \end{array}$$

*commutes and such that the primary obstruction to a section of $\pi$ is $\hat{c}_1$. Moreover the induced map between the fibres of the fibrations $B(\Pi, w) \to K(\Pi, 1)$ and $B^c(\Pi, w) \to K(\Pi, 1)$ is the usual fibration* BSpin $\to \mathrm{BSpin}^c$.

*Proof.* Consider the fibration $\pi \colon E \to B^c(\Pi, w)$ with fibre $K(\mathbb{Z}, 1)$ given by the class $\hat{c}_1 \in H^2(B^c(\Pi, w); \mathbb{Z})$ and the composition

$$p \circ \pi \colon E \longrightarrow B^c(\Pi, w) \longrightarrow \mathrm{BSO} \times K(\Pi, 1).$$

The restriction of $\pi$ to $\mathrm{BSpin}^c$ is the fibration given by $c_1$, which is the standard fibration BSpin $\to \mathrm{BSpin}^c$. To prove our claim it is therefore sufficient to show that $p \circ \pi$ and $q$ are fibre homotopy equivalent. Using the exact homotopy sequence for $\pi$ and its restriction to $\mathrm{BSpin}^c$ it is not difficult to see that $E$ has the same homotopy groups as $B(\Pi, w)$ and that the fibre $F$ of $p \circ \pi$ is a $K(\mathbb{Z}_2, 1)$. By construction, $(p \circ \pi)^*(w_2 + w) = 0$. As the kernel of the pullback is generated by the primary obstruction $\sigma$ to a section, this proves $\sigma = w_2 + w$ and our claim once we can show that $\sigma \neq 0$. But the existence of a section of $p \circ \pi$ would contradict $\pi_2(E) = 0$ and therefore $\sigma \neq 0$. □

**Remark 1.** The fibration $\pi$ is perhaps not uniquely determined by the class $\hat{c}_1$, but any two such fibrations differ only by a fibre homotopy equivalence over $B^c(\Pi, w)$.

Following the usual terminology, we will call a lift $\nu$ of the stable normal Gauss map $X \to BSO$ to $B^c(\Pi, w)$ a $B^c(\Pi, w)$–structure on a manifold $X$, similarly a lift to $B(\Pi, w)$ will be called a $B(\Pi, w)$–structure.

**Definition 7.** If $(\Pi, w, \hat{c}_1)$ is an admissible triple and $\nu$ is a $B^c(\Pi, w)$–structure on a 4–manifold $X$, we denote the class $v^*(\hat{c}_1) \in H^2(X; \mathbb{Z})$ by $c_1(\nu)$. This class will be referred to as the **Chern class** of $\nu$.

Clearly the explicit description of $\pi_2(B^c(\Pi, w))$ given in Lemma 4 shows that a $B^c(\Pi, w)$–structure $\nu$ on a manifold $X$ is a 2–equivalence if and only if the map

$$\nu_*: \pi_1(X) \longrightarrow \pi_1(B^c(\Pi, w)) = \Pi$$

is an isomorphism and there is a spherical homology class $\omega \in H_2(X; \mathbb{Z})$ such that $\langle c_1(\nu), \omega \rangle = 1$. Also note that a $B^c(\Pi, w)$–structure $\nu$ can be lifted to a $B(\Pi, w)$–structure if and only if $c_1(\nu) = 0$.

We can define a bordism group $\Omega_*^c(\Pi)$ of manifolds with $B^c(\Pi, w)$–structures in the usual way, the reader is referred to [10] for details. The inclusion of the fibre $B\mathrm{Spin}^c$ of the fibration $p'$ induces a homomorphism

$$\Omega_4^{\mathrm{Spin}^c} \longrightarrow \Omega_4^c(\Pi, w).$$

Using the fact that the 4–dimensional $\mathrm{Spin}^c$–bordism group is determined by the characteristic numbers $c_1^2$ and $\sigma$ [10], one easily verifies that this map is one–to–one. Also note that its image consists of exactly those bordism classes which admit simply connected representatives. It is sometimes convenient to work with the quotient

$$\widetilde{\Omega}_4^c(\Pi, w) = \Omega_4^c(\Pi, w) / \Omega_4^{\mathrm{Spin}^c}.$$

The homomorphism

$$ed: \widetilde{\Omega}_4^c(\Pi, w) \longrightarrow H_4(\Pi; \mathbb{Z})$$

mapping the bordism class of $(X, \nu)$ to $(p' \circ \nu)_*[X]$ will be called the edge homomorphism. Finally we denote the bordism group given by the fibration $q$ by $\Omega_*(\Pi, w)$. Again we have an inclusion $\Omega_4^{\mathrm{Spin}} \subset \Omega_4(\Pi, w)$. The fibration $\pi$ defines a homomorphism

$$\pi_*: \Omega_*(\Pi, w) \longrightarrow \Omega_*^c(\Pi)$$

mapping $\Omega_4^{\mathrm{Spin}}$ to $\Omega_4^{\mathrm{Spin}^c}$. Recall from Remark 1 that the fibration $\pi$ and therefore also the map $\pi_*$ may not be uniquely determined by $\hat{c}_1$, but a different choice of $\pi$ changes $\pi_*$ only by an automorphism of $\Omega_*(\Pi, w)$. In particular the image of $\pi_*$ is determined by $\hat{c}_1$, it consists of those classes which have a representative $(X, \nu)$ for which $\hat{c}_1(\nu) = 0$.

**Remark 2.** Given an admissible triple $(\Pi, w, \hat{c}_1)$, Lemma 5 shows that we have a commuting diagram

$$\begin{array}{ccccc}
\mathrm{BSpin} & \longrightarrow & B(\Pi, w) & \longrightarrow & K(\Pi, 1) \\
\downarrow & & \pi \downarrow & & id \downarrow \\
\mathrm{BSpin}^c & \longrightarrow & B^c(\Pi, w) & \longrightarrow & K(\Pi, 1)
\end{array}$$

of fibrations in which the vertical map on the left hand side is the usual fibration $\mathrm{BSpin} \to \mathrm{BSpin}^c$. As all the fibrations involved are orientable, we have, as explained in [11], James spectral sequences

$$E^2_{p,q} = H_p(\Pi; \Omega^{\mathrm{Spin}}_q) \Longrightarrow \Omega_{p+q}(\Pi, w)$$

and

$$\bar{E}^2_{p,q} = H_p(\Pi; \Omega^{\mathrm{Spin}^c}_q) \Longrightarrow \Omega^c_{p+q}(\Pi, w)$$

and a homomorphism between them, compatible with the map between the bordism groups induced by $\pi$

The relation between $B^c(\Pi, w)$–structures and embedded surfaces and the reason why we introduced the fibration $B^c(\Pi, w) \to \mathrm{BSO} \times K(\Pi, 1)$ is provided by the following observation.

**Lemma 6.** *Suppose that $X$ is a 4–manifold and that $\varphi \in H^2(X; \mathbb{Z})$ is a spherically characteristic cohomology class. Then there exists an admissible triple $(\Pi, w, \hat{c}_1)$ and a $B^c(\Pi, w)$–structure $\nu$ on $X$ such that $c_1(\nu) = PD(\varphi)$. Moreover we can assume that the homomorphism*

$$\nu_* \colon \pi_1(X) \longrightarrow \Pi$$

*is an isomorphism.*

*Proof.* By assumption, the pullback of $\varphi$ to the universal covering $\widetilde{X}$ of $X$ is characteristic, i.e. there is a cohomology class $v \in H^2(X; \mathbb{Z}_2)$ such that $\varphi \equiv w_2(X) + v$ and the pullback of $v$ is zero. Let $\Pi = \pi_1(X)$ and choose a classifying map $u \colon X \to K(\Pi, 1)$ for the universal covering of $X$. As the sequence

$$0 \longrightarrow H^2(\Pi; \mathbb{Z}_2) \xrightarrow{u^*} H^2(X; \mathbb{Z}_2) \longrightarrow H^2(\widetilde{X}; \mathbb{Z}_2)$$

is exact, there is a unique class $w \in H^2(\Pi; \mathbb{Z}_2)$ such that $v = u^* w$. Choose a class $\hat{c}_1 \in B^c(\Pi, w)$ such that $(\Pi, w, \hat{c}_1)$ is an admissible triple.

The product of the stable normal Gauss map of $X$ and $u$ defines a map $f \colon X \to \mathrm{BSO} \times K(\Pi, 1)$. The pullback $f^*(w_2 + w) = w_2(X) + v$ has a lift, namely $\varphi$, hence the Bockstein operator applied to that class is zero. This implies that we can find a lift $\nu_0 \colon X \to B^c(\Pi, w)$ of $f$, which is then a $B^c(\Pi, w)$–structure. Clearly

$$c_1(\nu_0) \equiv w_2(X) + u^* w \equiv \varphi \mod 2,$$

i.e. there is a class $d \in H^2(X; \mathbb{Z})$ such that $\varphi = c_1(\nu_0) - 2d$. By obstruction theory, there is another $B^c(\Pi, w)$–structure $\nu$ on $X$, whose projection to $\mathrm{BSO} \times K(\Pi, 1)$ is also $f$, such that the primary difference between these two lifts of $f$ is the class $d$. The exact homotopy sequence of the fibration $p$ shows that the map

$$\pi_2(K(\mathbb{Z}, 2)) = \mathbb{Z} \longrightarrow \pi_2(B^c(\Pi, w)) = \mathbb{Z}$$

induced by the inclusion of the fibre is the multiplication by 2, which, together with the explicit description of $\pi_2(B^c(\Pi, w))$ given in Lemma 4, implies that

$$c_1(\nu_0) - c_1(\nu) = 2d = c_1(\nu_0) - \varphi.$$

It follows that $c_1(\nu) = \varphi$. Moreover the composition $p' \circ \nu$ is $u$, and by the choice of $u$ the map $u_*: \pi_1(X) \to \Pi$ is an isomorphism, which implies that the same is true for the homomorphism induced by $\nu$. Hence the $B^c(\Pi, w)$–structure $\nu$ has all the desired properties and the proof is complete. $\square$

Roughly speaking, this lemma shows that the spherically characteristic classes can be characterized as the Chern classes of $B^c(\Pi, w)$–structures, although we have to use different values of $w$ to realize all spherically characteristic classes. In analogy with [11], one could call the pair $(\Pi, w)$ used in the proof the $w_2$-type of $\varphi$. We close this section by introducing a homomorphism which we will need later on.

**Definition 8.** Given an admissible triple $(\Pi, w, \hat{c}_1)$, we define a map

$$\hat{c}_1 \colon \widetilde{\Omega}_4^c(\Pi, w) \longrightarrow H_2(\Pi; \mathbb{Z})$$

by mapping a bordism class $[(X, \nu)]$ to $(p' \circ \nu)_* PD(c_1(\nu))$.

We will now see that this map can be partially computed, more precisely we can compute its restriction to the kernel of the edge homomorphism. As already mentioned above, there is, as in [11], a James spectral sequence

$$\bar{E}^2_{p,q} = H_p(\Pi; \Omega_q^{\text{Spin}^c}) \Longrightarrow \Omega^c_{p+q}(\Pi).$$

In this spectral sequence, the lowest term $F_0 \Omega_4^c(\Pi, w)$ in the filtration of $\Omega_4^c(\Pi, w)$ is $\Omega_4^{\text{Spin}^c}$. It is well known that $\Omega_3^{\text{Spin}^c} = 0 = \Omega_1^{\text{Spin}^c}$ and $\Omega_2^{\text{Spin}^c} = \mathbb{Z}$, which implies that we have a map

$$\bar{E}^2_{2,2} = H_2(\Pi; \mathbb{Z}) \longrightarrow \bar{E}^\infty_{2,2} \subset \widetilde{\Omega}_4^c(\Pi, w)$$

whose image is precisely the kernel of the edge homomorphism. To compute the composition of this map with $\hat{c}_1$, we will make use of the following naturality property, which will also be useful for other purposes.

**Lemma 7.** *Suppose that $\Pi$ and $\Pi'$ are groups and that we are given cohomology classes $w \in H^2(\Pi; \mathbb{Z}_2)$ and $w' \in H^2(\Pi; \mathbb{Z}_2)$. If $f \colon \Pi \to \Pi'$ is a homomorphism such that $f^* w' = w$, there is a map $\hat{f} \colon B^c(\Pi, w) \to B^c(\Pi', w')$ such that the diagram*

$$\begin{array}{ccc} B^c(\Pi, w) & \xrightarrow{\hat{f}} & B^c(\Pi', w') \\ \downarrow & & \downarrow \\ \text{BSO} \times K(\Pi, 1) & \xrightarrow{(id, f_*)} & \text{BSO} \times K(\Pi', 1) \end{array}$$

*commutes and such that the induced map $\text{BSpin}^c \to \text{BSpin}^c$ is homotopic to the identity over $\text{BSO}$. In particular $\hat{f}$ induces a map between the James spectral sequences for $\Omega_*^c(\Pi', w')$, and $\Omega_*^c(\Pi, w)$ such that the homomorphism between the $\bar{E}^2$–terms is the homomorphism induced by $f$.*

*Proof.* As, by construction, the map $(id, f_*)$ pulls back $w_2 + w'$ to $w_2 + w$, we can find a map

$$\hat{f} \colon B^c(\Pi, w) \longrightarrow B^c(\Pi', w')$$

such that the above diagram commutes. Choose classes $\hat{c}_1 \in H^2(B^c(\Pi, w); \mathbb{Z})$ and $\hat{c}_1' \in H^2(B^c(\Pi', w'))$ such that $(\Pi, w, \hat{c}_1)$ and $(\Pi', w', \hat{c}_1')$ are admissible triple. Now clearly the classes $\hat{f}^*\hat{c}_1'$ and $\hat{c}_1$ have the same reduction modulo two, namely $p^*(w_2 + w)$, and therefore we can, using the same argument as in the proof of Lemma 6, arrange for $\hat{f}^*\hat{c}_1' = \hat{c}_1$. This implies in particular that the induced self–map of BSpin$^c$ over BSO pulls back $c_1$ to itself, which, by obstruction theory, shows that it is vertically homotopic to the identity. □

**Lemma 8.** *Suppose that $(\Pi, w, \hat{c}_1)$ is an admissible triple. The composition*

$$\varphi \colon H_2(\Pi; \mathbb{Z}) = \bar{E}_{2,2}^2 \longrightarrow \bar{E}_{2,2}^\infty \subset \widetilde{\Omega}_4^c(\Pi, w) \xrightarrow{\hat{c}_1} H_2(\Pi; \mathbb{Z})$$

*is the multiplication by two.*

*Proof.* Suppose that we are given some $x \in H_2(\Pi; \mathbb{Z})$. Since we represent $x$ by a map from a surface of positive genus we can find an aspherical surface $F$ and a homomorphism $f \colon \pi_1(F) \to \Pi$ mapping the orientation to $x$. Let $w' = f^*w$. As $\beta(w')$ is zero, we can therefore use the naturality property stated in Lemma 7 to reduce to the case that $\beta(w) = 0$.

In this case, $B^c(\Pi, w) = $ BSpin$^c \times K(\Pi, 1)$, and $\hat{c}_1 = c_1 + \hat{w}$ for some class $\hat{w} \in H^2(\Pi; \mathbb{Z})$. Moreover a $B^c(\Pi, w)$–structure on a 4–manifold $X$ is the same as a pair $(s, u)$, where $s$ is a Spin$^c$–structure and $u$ is some map $X \to K(\Pi, 1)$, and the James spectral sequence is nothing but the Atiyah–Hirzebruch spectral sequence for the Spin$^c$–bordism group of $K(\Pi, 1)$. By the computations in [1], we only have to prove that the restriction of the map $\hat{c}_1$ to the kernel of the edge homomorphism is given by mapping a representative $(X, s, u)$ to $u_*PD(c_1(s))$, where $c_1(s) = s^*c_1$ is the Chern class of $s$. For that purpose suppose that $\alpha = [(X, s, u)]$ is in the kernel of $ed$, i.e. $u_*[X] = 0$. Then

$$\hat{c}_1(\alpha) = u_*PD(c_1(s) + u^*\hat{w}) = u_*PD(c_1(s)) + u_*([X] \cap u^*\hat{w})$$
$$= u_*PD(c_1(s)) + u_*[X] \cap \hat{w} = u_*PD(c_1(s))$$

which completes the proof. □

## 4. Stable embeddings of spherically characteristic surfaces

Having set up the necessary bordism theoretical framework in the preceding section, we can now start to prove our main results. As mentioned in the introduction, we will need the existence of certain 4–manifolds whose universal covering is spin, which we prove in the next lemma.

**Lemma 9.** *Suppose that $(\Pi, w, \hat{c}_1)$ is an admissible triple and that $w \neq 0$. Then there exists a $B^c(\Pi, w)$ 4–manifold $(X, \nu)$ such that $\sigma(X) \equiv 8 \mod 16$, the homomorphism*

$$\nu_* \colon \pi_1(X) \longrightarrow \Pi$$

*is an isomorphism and $c_1(\nu) = 0$. Moreover $(X, \nu)$ can be chosen such that its bordism class in $\widetilde{\Omega}_4^c(\Pi, w)$ is zero.*

*Proof.* The proof is a slight extension of the arguments in [11]. We will use the spectral sequences $E$ and $\bar{E}$ mentioned in Remark 2. It is well known (see for instance [10]) that

$$\Omega_1^{\text{Spin}} = \Omega_2^{\text{Spin}} = \mathbb{Z}_2 \ , \ \Omega_3^{\text{Spin}} = \Omega_3^{\text{Spin}^c} = \Omega_1^{\text{Spin}^c} = 0 \ , \ \Omega_2^{\text{Spin}^c} = \mathbb{Z}.$$

As all differentials emerging from $E^*_{2,2}$ and $\bar{E}^*_{2,2}$ are zero and $E^\infty_{1,3} = \bar{E}^\infty_{1,3} = 0$, we have a commuting diagram

$$\begin{array}{ccc} E^2_{2,2} = H_2(\Pi;\mathbb{Z}_2) & \longrightarrow & E^\infty_{2,2} \subset \Omega_4(\Pi,w)/\Omega_4^{\mathrm{Spin}} \\ \downarrow & & \downarrow \\ \bar{E}^2_{2,2} = H_2(\Pi;\mathbb{Z}) & \longrightarrow & \bar{E}^\infty_{2,2} \subset \Omega^c_4(\Pi,w)/\Omega_4^{\mathrm{Spin}^c} \end{array}$$

in which the vertical map on the left hand side is induced by $\Omega_2^{\mathrm{Spin}} \to \Omega_2^{\mathrm{Spin}^c}$ which is the trivial homomorphism, therefore this map itself is zero. As $w \neq 0$, the map $w_*\colon H_2(\Pi;\mathbb{Z}_2) \to H_2(K(\mathbb{Z}_2,2);\mathbb{Z}_2)$ is onto. Pick an element $x \in E^2_{2,2}$ such that $w_*(x) \neq 0$. Let

$$\alpha \in \Omega_4(\Pi,w)/\Omega_4^{\mathrm{Spin}}$$

denote its image under the upper horizontal map of the above diagram, and choose a representative $(X,\bar{\nu})$ of $\alpha$ such that $\bar{\nu}\colon X \to B(\Pi,w)$ is a 2–equivalence. As explained in [11], $X$ will then be a 4–manifold whose universal covering is spin, and $\sigma(X) \equiv 8 \mod 16$. The $B(\Pi,w)$–structure on $X$ induces a $B^c(\Pi,w)$–structure $(X,\nu)$ on $X$ such that $c_1(\nu) = 0$, and the fact that $\bar{\nu}$ is a 2–equivalence implies that the homomorphism $\nu_*\colon \pi_1(X) \to \Pi$ is an isomorphism, recall that $\pi$ induces an isomorphism between $\pi_1(B(\Pi,w))$ and $\pi_1(B^c(\Pi,w))$. Finally the commuting diagram above shows that the reduction of the bordism class

$$[(X,\nu)] \in \widetilde{\Omega}^c_4(\Pi,w) = \Omega^c_4(\Pi,w)/\Omega_4^{\mathrm{Spin}^c}$$

is zero and the proof is complete. □

**Lemma 10.** *Suppose that $(\Pi,w,\hat{c}_1)$ is an admissible triple and that $w \neq 0$. Then there is a $B^c(\Pi,w)$ 4–manifold $(X,\nu)$ representing $0 \in \Omega^c_4(\Pi,w)$ such that $\nu$ is a 2–equivalence and a $\pi_1$–null embedded surface $F \subset X$ representing $c_1(\nu)$ having Arf invariant one.*

*Proof.* Choose a $B^c(\Pi,w)$ 4–manifold $(X_0,\nu_0)$ as in the statement of Lemma 9. Then there is a simply connected $B^c(\Pi,w)$ 4–manifold $(X_1,\nu_1)$ such that

$$[(X_1,\nu_1)] = -[(X_0,\nu_0)] \in \Omega^c_4(\Pi,w).$$

We can also assume that $c_1(\nu_1)$ is a primitive class, if this is not the case we can add a copy of $\mathbb{C}P^2 \# \overline{\mathbb{C}P^2}$ with the unique $B^c(\Pi,w)$–structure having Chern class $(1,1)$. Let $X = X_0 \# X_1$. The $B^c(\Pi,w)$–structures $\nu_0$ and $\nu_1$ glue to give a $B^c(\Pi,w)$–structure $\nu$ on $X$, and as $X_1$ is simply connected, the map $\nu$ induces an isomorphism $\pi_1(X) = \Pi$. Note that $[(X,\nu)] = 0$ by the choice of $(X_1,\nu_1)$, in particular $c_1(\nu)^2 = 0 = c_1(\nu_1)^2$ and $\sigma(X) = 0$. Also observe that, as $c_1(\nu_1)$ is primitive and $X_1$ is simply connected, we can find a spherical class $\omega \in H_2(X;\mathbb{Z})$ such that $\langle c_1(\nu),\omega\rangle = 1$, consequently $\nu$ is a 2–equivalence.

Now pick a surface $F \subset X_1$ representing $c_1(\nu_1)$. Then this embedding defines a $\pi_1$–null embedding $F \hookrightarrow X$ representing $c_1(\nu)$ (recall that $c_1(\nu_0) = 0$), and as $F$ is entirely contained in the simply connected summand $X_1$ in which it is characteristic, we can compute its Arf invariant using the formula given in [3]. We obtain that

$$\mathrm{Arf}(F) \equiv \frac{1}{8}(c_1(\nu_1)^2 - \sigma(X_1)) \equiv \frac{1}{8}(\sigma(X) - \sigma(X_0)) \equiv 1 \mod 2$$

and the lemma is proved. □

*Proof of Theorem 6.* By Lemma 6, we can find an admissible triple $(\Pi, w, \hat{c}_1)$ and a $B^c(\Pi, w)$–structure $\nu$ on $X$ such that $PD(F) = c_1(\nu)$ and the induced homomorphism $\pi_1(X) \to \Pi$ is an isomorphism. Note that the existence of the class $\omega$ implies that $\nu$ is a 2–equivalence. By assumption, $c_1(\nu) = PD([F])$ is not characteristic, hence $w \neq 0$.

By Lemma 10, there is a $B^c(\Pi, w)$ 4–manifold $(X_0, \nu_0)$ representing $0 \in \Omega_4^c(\Pi, w)$ and a $\pi_1$–null embedding $F_0 \hookrightarrow X_0$ representing $c_1(\nu_0)$ such that $\text{Arf}(F_0) = 1$ and $\nu_0$ induces an isomorphism $\pi_1(X_0) = \Pi$. Let $X_1 = X \# X_0$ and glue the $B^c(\Pi, w)$–structures on $X$ and $X_0$ to obtain a $B^c(\Pi, w)$–structure $\nu_1$ on $X_1$. The connected sum of the surfaces $F$ and $F_0$ defines a $\pi_1$–null embedding $F_1 \hookrightarrow X_1$ representing $c_1(\nu_1)$ such that $\text{Arf}(F_1) = \text{Arf}(F) + 1$. The kernel of

$$(\nu_1)_* \colon \pi_1(X_1) = \pi_1(X) * \pi_1(X_0) \longrightarrow \Pi$$

is the subgroup which is normally generated by all elements $(\nu_*)^{-1}(g)((\nu_0)_*)^{-1}(g^{-1})$ where $g \in \Pi$. The group $\Pi$ is finitely generated, hence there are finitely many circles in $X_1 \setminus F_1$ whose homotopy classes normally generate the kernel of $(\nu_1)_*$. As $PD(F) = c_1(\nu_1)$, there exists a lift $\bar{\nu}$ of $\nu$ on $X_1 \setminus F_1$ to a $B(\Pi, w)$–structure. If we perform surgery along these circles compatible with that structure and equip the resulting manifold with the $B^c(\Pi, w)$–structure induced by an extension of $\bar{\nu}$, we obtain a $B^c(\Pi, w)$ 4–manifold $(X', \nu')$ representing the same bordism class as $(X_1, \nu_1)$ such that $\nu'$ induces an isomorphism $\pi_1(X') = \Pi$. As the surgery took place near circles disjoint from $F_1$, the surface $F_1$ survives to give a $\pi_1$–null embedded surface $F' \subset X'$ representing $c_1(\nu')$ which still has Arf invariant $\text{Arf}(F) + 1$. Furthermore we can represent the class $\omega \in H_2(X; \mathbb{Z}) \subset H_2(X_1; \mathbb{Z})$ by an immersed sphere disjoint from the $C_i$ which will also survive the surgery and gives an immersed sphere in $X'$ which has intersection number 1 with $F'$. Hence the $B^c(\Pi, w)$–structure $\nu'$ is a 2–equivalence.

Now we can apply the stable diffeomorphism classification due to M. Kreck. As the bordism class represented by $(X', \nu')$ is $[(X, \nu)]$, Theorem 2 in [6] shows that there are numbers $k, t$ and a diffeomorphism

$$\Phi \colon X \# k(S^2 \times S^2) \longrightarrow X' \# t(S^2 \times S^2)$$

which is compatible with the $B^c(\Pi, w)$–structures obtained by the $B^c(\Pi, w)$–structures on $X$ respectively $X'$ and the canonical $B^c(\Pi, w)$–structure on $S^2 \times S^2$. As the latter has Chern class zero, $\Phi^{-1}$ maps $F'$ to a $\pi_1$–null embedded surface in $X \# k(S^2 \times S^2)$ representing $c_1(\nu) = PD([F])$, and as the Arf invariant is clearly a diffeomorphism invariant, the Arf invariant of that surface will be $\text{Arf}(F) + 1$ as desired. □

As already mentioned in the introduction, we are now in the position to describe a bordism theoretical obstruction whose vanishing is a necessary and sufficient condition for a spherical and spherically characteristic homology class to be stably representable by an embedded sphere.

**Theorem 7.** *Suppose that $X$ is a 4–manifold and that $\xi \in H_2(X; \mathbb{Z})$ is a spherically characteristic spherical homology class. Choose an admissible triple $(\Pi, w, \hat{c}_1)$ and a $B^c(\Pi, w)$–structure $\nu$ on $X$ such that $c_1(\nu) = PD(\xi)$ and such that $\nu$ induces an isomorphism $\pi_1(X) = \Pi$. Consider the residue class*

$$\alpha(\nu) = [(X, \nu)] \in \widetilde{\Omega}_4^c(\Pi, w)/im(\Omega_4(\Pi, w) \longrightarrow \widetilde{\Omega}_4^c(\Pi, w))$$

of the bordism class given by $(X, \nu)$. Then the class $\xi$ can be stably represented by a $\pi_1$–null embedding if and only if $\alpha(\nu) = 0$.

If, in addition, $\xi$ is not characteristic and there is a spherical class $\omega$ such that $\xi \cdot \omega = 1$, then the class $\xi$ can be stably represented by an embedded sphere if and only if $\alpha(\nu) = 0$.

*Proof.* Theorem 6 implies that we only have to prove the first part of the statement. So let us suppose that the obstruction $\alpha(\nu)$ vanishes. The idea of the proof is to decompose the manifold $X$ stably into a simply connected manifold and a 4–manifold whose universal covering is spin such that the class $\xi$ is living in the simply connected part, then we can find a $\pi_1$–null embedded representative for trivial reasons. Again the main tool used for that purpose will be Kreck's modified surgery theory.

First let us assume that there is a spherical class $\omega$ such that $\omega \cdot \xi = 1$. Then the $B^c(\Pi, w)$–structure $\nu$ is a 2–equivalence. As $\alpha(\nu) = 0$, the reduced bordism class $[(X, \nu)] \in \widetilde{\Omega}_4^c(\Pi, w)$ is contained in the image of $\Omega_4(\Pi, w)$. Therefore we can find a simply connected $B^c(\Pi, w)$ 4–manifold $(X_1, \nu_1)$ and a $B^c(\Pi, w)$ 4–manifold $(X_0, \nu_0)$ such that

$$[(X, \nu)] = [(X_1, \nu_1)] + [(X_0, \nu_0)] \in \Omega_4^c(\Pi, w)$$

and $c_1(\nu_0) = 0$. We can also arrange for $(\nu_0)_* \colon \pi_1(X_0) \to \Pi$ to be an isomorphism. Furthermore we can assume that $c_1(\nu_1)$ is a primitive class, if this is not the case simply add a copy of $\mathbb{C}P^2 \# \overline{\mathbb{C}P^2}$ with the unique $B^c(\Pi, w)$–structure having Chern class $(1, 1)$ (note that this does not change the bordism class). Consider the 4–manifold $X' = X_0 \# X_1$ together with the $B^c(\Pi, w)$–structure $\nu'$ defined by $\nu_0$ and $\nu_1$. Then $c_1(\nu') = c_1(\nu_1)$ has support in the simply connected summand $X_0$ and can therefore be represented by a $\pi_1$–null embedded surface $F' \subset X'$. Moreover the $B^c(\Pi, w)$–structure $\nu'$ is a 2–equivalence. Applying once more Theorem 2 in [6] now shows that there is a diffeomorphism

$$\Phi \colon X \# s(S^2 \times S^2) \longrightarrow X' \# t(S^2 \times S^2)$$

for certain non–negative numbers $s$ and $t$ which is compatible with the $B^c(\Pi, w)$–structures and therefore maps $\xi$ to $PD(c_1(\nu'))$. The preimage $F$ of $F'$ is then the desired representative of $\xi$. We remark that this embedding has the additional property of being $\pi_1$–negligible, i.e. $\pi_1(X_k \setminus F) = \pi_1(X_k)$.

As to the general case, note that for every spherical and spherically characteristic class $\xi \in H_2(X; \mathbb{Z})$, the class

$$(\xi, 1) \in H_2(X \# \mathbb{C}P^2; \mathbb{Z}) = H_2(X; \mathbb{Z}) \oplus \mathbb{Z}$$

is again spherical and spherically characteristic. Moreover there is a spherical class $\omega$ such that $\omega \cdot (\xi, 1) = 1$, we can choose $\omega = (0, 1)$. We can also extend $\nu$ to a $B^c(\Pi, w)$–structure $\nu'$ on $X \# \mathbb{C}P^2$ such that $\hat{c}_1(\nu') = PD((\xi, 1))$ and $\alpha(\nu') = \alpha(\nu) = 0$. Therefore we can represent the class $(\xi, 1)$ stably by a $\pi_1$–null embedded surface and then apply the results in [1] to conclude that the same is true for $\xi$.

Now let us assume that conversely the class $\xi$ can be represented by a $\pi_1$–null embedded surface $F \subset X_k$ for some $k$. As we can extend the $B^c(\Pi, w)$–structure on $X$ to $X_k$ by the canonical $B^c(\Pi, w)$–structure on $S^2 \times S^2$ without changing the bordism class, we can assume that $k = 0$, i.e. $X = X_k$. If the Arf invariant of $F$ is zero we can, after adding additional copies of $S^2 \times S^2$, also assume that $F$ is a sphere. By taking connected sum with an embedded sphere in an appropriate simply connected $B^c(\Pi, w)$ 4–manifold, we can finally reduce to the case that $c_1(\nu)^2 = F \cdot F = 1$. Now the idea is to blow down the sphere $F$ to show that $[(X, \nu)]$ is in the image of $\Omega_4(\Pi, w)$ and therefore $\alpha(\nu) = 0$. More precisely, $c_1(\nu)$ is Poincaré dual to $[F]$, hence it vanishes on the complement of a tubular neighborhood $T$ of $F$. Thus we can lift $\nu$ to a $B(\Pi, w)$–structure on $X \setminus T$. As the sphere $F$ has self–intersection number one, we can blow down $F$ to obtain a decomposition

$$X = Y \# \mathbb{C}P^2$$

for some 4–manifold $Y$, i.e. $Y$ is the manifold which is obtained from $X$ by deleting $T$ and gluing in a 4–dimensional ball instead. Furthermore we can find $B^c(\Pi, w)$–structures $\nu_0$ on $Y$ and $\nu_1$ on $\mathbb{C}P^2$ such that

$$[(X, \nu)] = [(Y, \nu_0)] + [(\mathbb{C}P^2, \nu_1)] \in \widetilde{\Omega}_4^c(\Pi, w).$$

Since $\mathbb{C}P^2$ is simply connected, the bordism class of $(\mathbb{C}P^2, \nu_1)$ in $\widetilde{\Omega}_4^c(\Pi, w)$ is zero. Clearly a lift of $\nu$ to a $B(\Pi, w)$–structure on $X \setminus T$ can be extended to a lift of $\nu_0$, and therefore we obtain that $[(X, \nu)]$ is in the image of the homomorphism

$$\pi_* \colon \Omega_4(\Pi, w) \longrightarrow \widetilde{\Omega}_4^c(\Pi, w)$$

induced by $\pi$, i.e. $\alpha(\nu) = 0$.

If the Arf invariant of $F$ is one, we can take the connected sum with $\mathbb{C}P^2$ and the $B^c(\Pi, w)$–structure whose first Chern class is three times the generator to obtain another representative $(X \# \mathbb{C}P^2, \nu')$ of $\alpha(\nu)$. As $PD(c_1(\nu')) = [F] + 3\gamma$, where $\gamma \in H_2(\mathbb{C}P^2; \mathbb{Z})$ denotes the generator, and $3\gamma$ can be represented by an embedded surface having Arf invariant 1, we can glue that surface with $F$ to obtain another $\pi_1$–null embedded surface $F'$ representing $PD(c_1(\nu'))$ which has Arf invariant zero. Then we can argue as above to conclude that $\alpha(\nu)$ vanishes. □

Note that not every element in the quotient $\widetilde{\Omega}_4^c(\Pi, w)/im(\Omega_4(\Pi, w))$ can actually be realized as an obstruction $\alpha(\nu)$. The reason is that for a $B^c(\Pi, w)$–structure $\nu$ on a 4–manifold $X$ such that $\nu_* \colon \pi_1(X) \to \Pi$ is an isomorphism, $PD(c_1(\nu))$ is spherical if and only if $\hat{c}_1([(X, \nu)]) = 0$, this follows immediately from Definition 8. Using that it is clear that the residue classes which can be realized as obstructions $\alpha(\nu)$ are exactly those which are in the kernel of $\hat{c}_1$, i.e. our obstructions live in the group

$$H(\Pi, w, \hat{c}_1) = \ker(\hat{c}_1)/im(\Omega_4(\Pi, w) \longrightarrow \widetilde{\Omega}_4^c(\Pi, w))$$

which measures the extent to which the sequence

$$\Omega_4(\Pi, w) \longrightarrow \widetilde{\Omega}_4^c(\Pi, w) \xrightarrow{\hat{c}_1} H_2(\Pi; \mathbb{Z})$$

fails to be exact. The main result of [1] is that this group is zero if $w = 0$. To prove Theorem 4, we only have to exhibit an admissible triple $(\Pi, w, \hat{c}_1)$ for which $H(\Pi, w, \hat{c}_1)$ is not zero, then our claim follows from Theorem 7. For that purpose we need two technical lemmas.

**Lemma 11.** *For every $w \in H^2(\mathbb{Z} \oplus \mathbb{Z}_2; \mathbb{Z}_2)$ and for every choice of $\hat{c}_1$ such that $(\mathbb{Z} \oplus \mathbb{Z}_2, w, \hat{c}_1)$ is an admissible triple, there is a non–zero element $\alpha$ in the kernel of the edge homomorphism*

$$ed \colon \widetilde{\Omega}_4^c(\mathbb{Z} \oplus \mathbb{Z}_2, w) \longrightarrow H_4(\mathbb{Z} \oplus \mathbb{Z}_2; \mathbb{Z})$$

*such that $\hat{c}_1(\alpha) = 0$.*

*Proof.* We will again use the spectral sequence
$$\bar{E}^2_{p,q} = H_p(\mathbb{Z} \oplus \mathbb{Z}_2; \Omega_q^{\mathrm{Spin}^c}) \Longrightarrow \Omega_{p+q}^c(\mathbb{Z} \oplus \mathbb{Z}_2, w).$$
As $\bar{E}^2_{1,3} = \bar{E}^2_{3,1} = 0$, the kernel of the edge homomorphism on $\widetilde{\Omega}_4^c(\mathbb{Z} \oplus \mathbb{Z}_2, w)$ is $\bar{E}^\infty_{2,2}$, which is the quotient of $\bar{E}^2_{2,2}$ by the image of the differential $d_3$. Let $\iota\colon \mathbb{Z}_2 \to \mathbb{Z} \oplus \mathbb{Z}_2$ denote the natural inclusion. Using again the naturality property expressed in Lemma 7, we can compare $\bar{E}$ with the corresponding spectral sequence for $\Omega_*^c(\mathbb{Z}_2, \iota^*w)$ to see that this differential is zero. Consequently
$$\bar{E}^\infty_{2,2} = \bar{E}^2_{2,2} = H_2(\mathbb{Z} \oplus \mathbb{Z}_2; \mathbb{Z}) = \mathbb{Z}_2.$$
Now let $x$ denote the non–zero element of that group and let $\alpha$ denote its image under the map
$$\bar{E}^2_{2,2} \longrightarrow \bar{E}^\infty_{2,2} \subset \widetilde{\Omega}_4^c(\mathbb{Z} \oplus \mathbb{Z}_2, w).$$
Then $\alpha$ is a non–zero element in the kernel of $ed$. By Lemma 8, the composition
$$H_2(\mathbb{Z} \oplus \mathbb{Z}_2; \mathbb{Z}) = \bar{E}^2_{2,2} \longrightarrow \widetilde{\Omega}_4^c(\mathbb{Z} \oplus \mathbb{Z}_2, w) \xrightarrow{\hat{c}_1} H_2(\mathbb{Z} \oplus \mathbb{Z}_2; \mathbb{Z})$$
is multiplication by 2. Thus $\hat{c}_1(\alpha) = 2x = 0$ and the lemma is proved. $\square$

**Lemma 12.** *Let $t \in H^1(\mathbb{Z} \oplus \mathbb{Z}_2; \mathbb{Z}_2)$ denote the pullback of the non–zero element in $H^1(\mathbb{Z}_2; \mathbb{Z}_2)$. For every choice of $\hat{c}_1$ such that $(\mathbb{Z} \oplus \mathbb{Z}_2, t^2, \hat{c}_1)$ is an admissible triple, the map*
$$\Omega_4(\mathbb{Z} \oplus \mathbb{Z}_2, t^2) \longrightarrow \widetilde{\Omega}_4^c(\mathbb{Z} \oplus \mathbb{Z}_2, t^2)$$
*maps the kernel of the edge homomorphism to zero.*

*Proof.* We will use the abbreviations $\Pi = \mathbb{Z} \oplus \mathbb{Z}_2$ and $w = t^2$. Let us again consider the spectral sequences $E$ and $\bar{E}$ and the homomorphism $E \to \bar{E}$ as described in Remark 2. We have to prove that the subgroup $F_3\Omega_4(\Pi, w)$ of $\Omega_4(\Pi, w)$, which is exactly the kernel of the edge homomorphism, is mapped to the subgroup
$$\Omega_4^{\mathrm{Spin}^c} = F_0\Omega_4^c(\Pi, w) \subset \Omega_4^c(\Pi, w).$$
Clearly $\bar{E}^\infty_{3,1} = \bar{E}^\infty_{1,3} = 0$, so
$$F_3\Omega_4^c(\Pi, w)/\Omega_4^{\mathrm{Spin}^c} = \bar{E}^\infty_{2,2}.$$
By the computations in [11], the differential
$$d_2\colon E^2_{3,1} = H_3(\Pi; \mathbb{Z}_2) \longrightarrow E^2_{1,2} = H_1(\Pi; \mathbb{Z}_2)$$
is dual to the cup product with $w = t^2$ which is easily seen to be an isomorphism. Hence $E^\infty_{3,1} = 0$. Obviously also $E^\infty_{1,3}$ vanishes, so we have
$$F_3\Omega_4(\Pi, w)/\Omega_4^{\mathrm{Spin}} = E^\infty_{2,2}.$$
Moreover all differentials emerging from $E^*_{2,2}$ and $\bar{E}^*_{2,2}$ are zero and therefore we have the commuting diagram
$$\begin{array}{ccc} \bar{E}^2_{2,2} & \longrightarrow & F_3\Omega_4^c(\Pi, w)/\Omega_4^{\mathrm{Spin}^c} = \bar{E}^\infty_{2,2} \\ \uparrow & & \uparrow \\ E^2_{2,2} & \longrightarrow & F_3\Omega_4(\Pi, w)/\Omega_4^{\mathrm{Spin}} = E^\infty_{2,2} \end{array}$$
in which both horizontal maps are surjections. But as already observed in the proof of Lemma 9, the vertical map on the left hand side is zero, and therefore also the vertical map on the right hand side is zero. This completes the proof of the lemma. □

*Proof of Theorem 4.* Let $\Pi = \mathbb{Z} \oplus \mathbb{Z}_2$ and $w = t^2$ as in Lemma 12. Choose a class $\hat{c}_1$ such that $(\Pi, w, \hat{c}_1)$ is an admissible triple. By Lemma 11, there exists a non–zero element $\alpha$ in the kernel of the edge homomorphism which is also in the kernel of $\hat{c}_1$. Now $\alpha$ is not in the image of the map

$$\pi_* \colon \Omega_4(\Pi, w) \longrightarrow \widetilde{\Omega}_4^c(\Pi, w),$$

because every preimage of $\alpha$ had to be in the kernel of the edge homomorphism which is mapped to zero by Lemma 12. Therefore $\alpha$ is a non–zero element in the group

$$H(\Pi, w, \hat{c}_1) = \ker(\hat{c}_1)/im(\Omega_4(\Pi, w) \longrightarrow \widetilde{\Omega}_4^c(\Pi, w)).$$

Choose a representative $(X, \nu)$ of $\alpha$ such that $\nu$ is a 2–equivalence. Then the class $c_1(\nu)$ is spherically characteristic and $\pi_1(X) = \Pi$. Moreover $(p' \circ \nu)_* \xi = \hat{c}_1(\alpha) = 0$. Since the map $u = p' \circ \nu$ is a 2–equivalence, the sequence ("Hopf–sequence")

$$\pi_2(X) \longrightarrow H_2(X; \mathbb{Z}) \xrightarrow{u_*} H_2(\Pi; \mathbb{Z}) \longrightarrow 0$$

is exact, which implies that $\xi$ is spherical. By Theorem 7, the fact that $\alpha = \alpha(\nu)$ is not zero shows that the class $\xi$ cannot be stably represented by a $\pi_1$–null embedding, and the proof is complete. □

Although we have just seen that there may be obstructions to representing homology classes stably by $\pi_1$–null embeddings, there are conditions on the fundamental group which guarantee that these obstructions vanish.

**Theorem 8.** *Suppose that $X$ is a 4–manifold such that $H_2(\pi_1(X); \mathbb{Z})$ does not contain 2–torsion and $H_4(\pi_1(X); \mathbb{Z})$ is finite of odd order. Then every spherical class $\xi \in H_2(X; \mathbb{Z})$ can be stably represented by a $\pi_1$–null embedding.*

Examples of groups which fulfill the hypothesis of Theorem 8 are free groups, cyclic groups, the fundamental groups of surfaces and the fundamental groups of aspherical 3–manifolds. We also mention that there are further groups for which explicit computation shows that the obstructions $\alpha(\nu)$ are always zero, but which do not fulfill the conditions of the theorem, for instance $(\mathbb{Z}_2)^r$.

*Proof.* By Theorem 3 and Theorem 7, we only have to prove that given an admissible triple $(\Pi, w, \hat{c}_1)$ for which $H_2(\Pi; \mathbb{Z})$ does not contain 2–torsion and $H_4(\Pi; \mathbb{Z})$ is finite of odd order, the kernel of $\hat{c}_1$ coincides with the image of $\Omega_4(\Pi, w)$. Using the James spectral sequence for $\Omega_*^c(\Pi, w)$, one can easily derive an exact sequence

$$H_2(\Pi; \mathbb{Z}) \longrightarrow \widetilde{\Omega}_4^c(\Pi, w) \xrightarrow{ed} H_4(\Pi; \mathbb{Z})$$

in which the map on the right hand side is the usual edge homomorphism. Now assume that we are given an element $\alpha \in \widetilde{\Omega}_4^c(\Pi, w)$ which is in the kernel of $\hat{c}_1$. In the James spectral sequence for $\Omega(\Pi, w)$, all the differentials emerging from $E_{4,0}^*$ are zero, because they map subgroups of $H_4(\Pi; \mathbb{Z})$ to finite abelian 2–groups and we assume $H_4(\Pi; \mathbb{Z})$ to be finite of odd order. Therefore the edge homomorphism

$$ed \colon \Omega_4(\Pi, w) \to H_4(\Pi; \mathbb{Z})$$

is onto. By subtracting an element $\beta$ in the image of $\Omega_4(\Pi, w)$ such that $ed(\alpha) = ed(\beta)$, we can therefore reduce to the case that the edge homomorphism maps $\alpha$ to zero. By the exactness of the sequence above, this implies that $\alpha$ is in the image of $H_2(\Pi;\mathbb{Z}) \to \widetilde{\Omega}_4^c(\Pi,w)$. Pick a preimage $x \in H_2(\Pi;\mathbb{Z})$. Then we have

$$2x = \varphi(x) = \hat{c}_1(\alpha) = 0,$$

where $\varphi$ is as in Lemma 8. By assumption, $H_2(\Pi;\mathbb{Z})$ does not have 2–torsion, and therefore we can conclude that $x = 0$ which proves the Theorem. $\square$

Mathematisches Institut, Theresienstr. 39, 80333 München, Germany
*E-mail address*: bohr@mathematik.uni-muenchen.de